\documentclass[fleqn]{mat01}
\usepackage{times,mathtimy,amssymb}
\begin{document}

\setcounter{page}{53}
\firstpage{53}

\font\xx=msam5 at 9pt
\def\ab{\mbox{\xx{\char'03}}}

\font\sa=tibi at 10.4pt

%\setcounter{page}{1}
%\firstpage{543}

\def\defi{\trivlist\item[\hskip\labelsep{\bf DEFINITION.}]}
\def\remark{\trivlist\item[\hskip\labelsep{\it Remark.}]}
\def\remarks{\trivlist\item[\hskip\labelsep{\it Remarks}]}
\def\noot{\trivlist\item[\hskip\labelsep{{\it Note.}}]}
\def\thoe{\trivlist\item[\hskip\labelsep{{\bf Theorem}}]}
\newtheorem{theo}{Theorem}
\renewcommand\thetheo{\arabic{section}.\arabic{theo}}
\newtheorem{theor}{\bf Theorem}
\newtheorem{lem}{Lemma}
\newtheorem{propo}{\rm PROPOSITION}
\newtheorem{rema}{Remark}
\newtheorem{exam}{Example}
\newtheorem{coro}{\rm COROLLARY}

\title{Very smooth points of spaces of operators}

\markboth{T S S R K Rao}{Very smooth points of spaces of operators}

\author{T S S R K RAO}

\address{Indian Statistical Institute, R.V. College Post, Bangalore 560 059, India\\
\noindent E-mail: tss@isibang.ac.in}

\volume{113}

\mon{February}

\parts{1}

\Date{MS received 18 August 2001; revised 22 January 2002}

\begin{abstract}
In this paper we study very smooth points of Banach spaces with special
emphasis on spaces of operators. We show that when the space of compact
operators is an $M$-ideal in the space of bounded operators, a very
smooth operator $T$ attains its norm at a unique vector $x$ (up to a
constant multiple) and $T(x)$ is a very smooth point of the range space.
We show that if for every equivalent norm on a Banach space, the dual
unit ball has a very smooth point then the space has the
Radon--Nikod\'{y}m property. We give an example of a smooth Banach space
without any very smooth points. 
\end{abstract}


\maketitle

\section{Introduction} 

A Banach space $X$ is said to be {\it very smooth} if every unit vector
has a unique norming element in $X^{***}$ (here $X$ is being
considered as a subspace of $X^{**}$ under the canonical
embedding, see [S]). In this paper we study a local version of the
notion of `very smooth space' by calling a unit vector of $X$ a very
smooth point if it is also a smooth point of $X^{**}$ (recall that
a unit vector is a smooth point if it has a unique norming functional in
the dual). These notions are related to differentiability of the norm at
these points, see [S]. In particular for the space of compact operators
we will be considering differentiability in the direction of every
bounded operator.

Identification of smooth points of spaces of operators ,
$C^{*}$-algebras and their generalizations has received a lot of
attention in the literature. See [H, GY1, GY2, KY, MR, TW,
W1, W2, HWW] Chapter VI. Motivated by these results, in this
paper we study very smooth points of spaces of operators and
$JB^{*}$-triples. 

To do this we first prove a proposition involving the notion of an
$M$-ideal that allows us to `lift' very smooth points from subspaces. We
use this to show that depending on the `position' of an $M$-ideal, all
smooth points of a Banach space can be very smooth. We give several
examples from spaces of operators where these ideas apply. 

In order to study very smooth points of the space of compact
operators ${\cal K}(X,Y)$ we make use of a characterization due to
Heinrich [H] of smooth points of this space (see \S2). Let $X_1$
denote the closed unit ball of $X$. We recall from [S] and [GI] that a
unit vector $x \in X_1$ is a very smooth point if and only if its unique
norming functional $x^{*}$ is a point of continuity for the identity
map $i : (X^{*}_1, w^{*}) \rightarrow (X^{*}_1, weak)$. Our approach
involves studying these points of $w^{*}-w$ continuity. For 
${\cal K}(X,Y)$ we show that a very smooth point attains its norm and the image
vector is a very smooth point. When $Y = C(K)$ for a compact space $K$
we could give a complete description of very smooth points. More
generally we show that $f \in C(K,X)$ (the space of $X$-valued
continuous functions on $K$, equipped with the supremum norm) is a very
smooth point iff there is a unique isolated point $k$ such that $f(k)$
is a very smooth point of $X$. We also show that any very smooth point
of a $JB^{*}$-triple is a Fr\'{e}chet smooth point. We give a
description of very smooth points of ${\cal L}(X,Y)$ when ${\cal
K}(X,Y)$ is an $M$-ideal in ${\cal L}(X,Y)$. See Chapter VI of [HWW] and
[KW] for several examples of Banach spaces $X$ and $Y$ for which ${\cal
K}(X,Y)$ is an $M$-ideal in ${\cal L}(X,Y)$. 

In the third section of this paper we prove a result analogous to a
result of Ruess and Stegall ([RS1, RS2]) by showing that when ${\cal K}(X,Y)$ is
an $M$-ideal in ${\cal L}(X,Y)$, then the dual unit balls of both the
spaces have the same points of $w^{*}-w$ continuity. 

In the final section of the paper we consider spaces that fail to
have very smooth points. Considering the inclusion $X \subset Y \subset
X^{**}$ under the canonical embedding, we show that if $X$ is a
proper $L$-ideal in $Y$ then $X_1$ has no very smooth points. This
allows us to give an example of a smooth space in which no unit vector
is very smooth. We complement our earlier work from [R2] and [R3] by
showing that for an infinite compact set $K$ and for a reflexive Banach
space $X$, if ${\cal L}(X,C(K))$ is a dual space then there is no point
of $w^{*}-w$ continuity in the unit ball. 
 
Our notation and terminology is standard and can be found in [HWW]. For
a Banach space $X$ by $\partial_e X_1$ we denote the set of extreme
points.

\section{Very smooth points} 

Let $M \subset X$ be a closed subspace. It was observed in [MR]
that if $x \in M$ is a smooth point of $X$ then it is a smooth point of
$M$. It is easy to see that if every continuous linear functional on $M$
has a unique norm preserving extension to $X$ then every smooth point of
$M$ is also a smooth point of $X$. Our first result addresses this
question for very smooth points and involves the notion of an $M$-ideal.
We recall from Chapter I of [HWW] that a closed subspace $M \subset X$
is an $M$-ideal, if there is a projection $P$ in ${\cal L}(X^{*})$ such
that $\|P(x^{*})\| + \|x^{*} - P (x^{*})\| = \|x^{*}\|$ for all
$x^{*} \in X^{*}$ and $\hbox{ker}\ P = M^\bot$. We note from 
Proposition~1.1.12 in [HWW] that if $M$ is an $M$-ideal in $X$ then continuous
linear functionals on $M$ have unique norm preserving extensions to $X$.
We also note that in $C^{*}$-algebras $M$-ideals are precisely closed
two sided ideals. We first prove a lemma on smooth points in the
$\ell^{\infty}$-direct sum of two spaces.

\begin{lem}
Suppose $X = M \oplus_{\infty} N $. Let $m \in M, \|m\| = 1$ and $n \in
N, \|n\| < 1$. $x = m+n$ is a smooth point of $X$ if and only if $m$ is
a smooth point of $M$. 
\end{lem}

\begin{proof}
Suppose $x$ is smooth in $X$. Let $\|m^{*}_i\|=1=m^{*}_i(m) $,
$m^{*}_i \in M^{*}$ for $i=1,2$. Then $(m^{*}_1, 0), (m^{*}_2, 0)$
are two linear functionals in $X^{*} = M^{*} \oplus_1 N^{*}$
($\ell^1$-direct sum) that attain their norm at $x$. Hence
$m^{*}_1=m^{*}_2$ and therefore $m$ is smooth in $M$.

Conversely suppose $m$ is smooth in $M$. Suppose
$\|x_i^{*}\|=1=x_i^{*} (x)$ for $i=1,2$. Since $X^{*} = M^\bot
\oplus_1 N^\bot = M^{*} \oplus_1 N^{*}$, $x_i^{*}=m_i^{*}+n_i^{*}$
and $1= \|m_i^{*} \| + \|n_i^{*}\| = m_i^{*} (m)+n_i^{*} (n)$.
Suppose $n_i^{*} \neq 0$. Thus $1= \|m_i^{*}\|
{m_i^{*}(m)}/{\|m_i^{*}\|} + \|n_i^{*}\|
{n_i^{*}(n)}/{\|n_i^{*}\|}$ implies $\|n\|=1$. Hence $n_i^{*}=0$
for $i=1,2$ and $m_1^{*}=m_2^{*}$. Therefore $x$ is smooth in $X$.
\hfill \ab
\end{proof}

\setcounter{propo}{1}
\begin{propo}$\left.\right.$\vspace{.3pc}

\noindent Let $m \in M$ be a very smooth point of $X$ then it
is a very smooth point of $M$. If $M \subset X$ is an $M$-ideal then a
very smooth point of $M$ is also a very smooth point of $X$.
\end{propo}

\begin{proof}
Let $m \in M$ be a very smooth point of $X$ . Since
$M^{{*}{*}} = M^{\bot\bot} \subset X^{{*}{*}}$ under the canonical
embedding and since $m$ is a smooth point of $X^{{*}{*}}$, by applying
Lemma 2.1 of [MR] we get that $m$ is a smooth point of $M^{{*}{*}}$.
Thus $m$ is a very smooth point of $M$.

Assume further that $M$ is an $M$-ideal. We have from the definition
that $X^{{*}{*}} = M^{{*}{*}} \oplus_{\infty} (M^{*})^{\bot}$. Now
if $m \in M$ is a very smooth point then it is a smooth point of
$M^{{*}{*}}$ and hence it follows from our lemma that it is a smooth
point of $X^{{*}{*}}$ and hence a very smooth point of $X$.  \hfill \ab
\end{proof}

\setcounter{rema}{2}
\begin{rema}
{\rm Part of Proposition~2 is an abstract version of Proposition 3.2
in [GI], where the authors showed by different methods that if $X$ is an
$M$-ideal in its bidual (under the canonical embedding), then very
smooth points of $X$ are also very smooth in $X^{{*}{*}}$. The
$M$-ideal condition in the above proposition cannot be replaced by `$M$
is an ideal in $X$' (we recall from [GKS] that $M \subset X$ is an ideal
if there is a projection $P$ of norm one in the dual such that $\hbox{ker} P =
M^{\bot}$ ). In [GI] the authors have constructed an example of a Banach
space $X$ and a very smooth point of $X$ that is not very smooth in
$X^{{*} {*}}$ ($X$ is canonically embedded in $X^{{*} {*}}$ and
hence an ideal). 
}
\end{rema}

However, for $M \subset X$ if one assumes that there is a norm-one
projection $P:X^{*} \rightarrow X^{*}$ such that ker $P=M^\perp$ and
$P(X^{*})_1$ is $w^{*}$-dense in $X^{*}_1$ then since under these
conditions $M \subset X \subset M^{{*} {*}}$ (under the canonical
embedding, see [R6]), we have that a very smooth point of $M$ is a
smooth point of $X$. 

If $X^{*}$ or $Y$ has the compact metric approximation property
(CMAP), then it is easy to see that $K(X,Y) \subset {\cal L} (X,Y)$
satisfies the above condition (see [R6]). Thus for a compact operator
that is a very smooth point, directional derivatives exist in the
direction of all bounded operators.

Smooth points of operator spaces has been extensively studied in the
literature. We now undertake to study very smooth points of these
spaces. As noted before one of the motivations is that under the
assumption of the CMAP, a very smooth point of $K(X,Y)$ is a smooth
point of ${\cal L}(X,Y)$. 

Before doing this we indicate one more proposition that exhibits the
presence of very smooth points depending on the position of an $M$-ideal
and involves the notion of a Hahn--Banach smooth space considered in [S].

To study very smooth points in this set up we use a characterization
of very smooth points obtained in [GI], as smooth points for which the
unique norming linear functional in $X^{*}_1$ is also a point of
continuity for the identity map on $X^{*}_1$ equipped with the
weak$^{*}$ and weak topologies on the domain and the range respectively. 

\setcounter{propo}{3}
\begin{propo}$\left.\right.$\vspace{.3pc}

\noindent Let $J$ be a Hahn--Banach smooth space. Suppose $J \subset X$ is an
$M$-ideal and $X/J$ does not have any smooth points. Then every smooth
point of $X$ is a very smooth point. 
\end{propo}

\begin{proof}
Let $x \in X$ be a smooth point. Since $X/J$ does not have
any smooth points, and since $J$ is an $M$-ideal, there exists a unique
$j^{*} \in J^{*}_1$ such that $j^{*} (x)=1$. Since $J$ is Hahn--Banach
smooth, $j^{*}$ is a point of $w^{*} - w$ continuity of $J_1^{*}$
(see [HWW], Lemma~III.2.14). We now claim that $j^{*}$ is a point of
$w^{*} -w$ continuity of $X_1^{*}$. Let $\{x_{\alpha}^{*}\}_{\alpha
\in I} \subset X^{*}_1$ be a net converging to $j^{*}$ in the
weak$^{*}$-topology of $X^{*}$. Clearly the net
$\{P(x^{*}_{\alpha})\}_{\alpha \in I}$ converges to $j^{*}$ in the
weak$^{*}$-topology of $J^{*}$. In view of our assumption on $J$,
this convergence is thus in the weak topology. Also $ 1 = \lim
\|P(x^{*}_{\alpha})\| = \lim \|x^{*}_{\alpha}\|$. By the defining
property of $P$ we get that the net $\{\|x^{*}_{\alpha} -
P(x^{*}_{\alpha})\|\}_{\alpha \in I} $ converges to $0$. Now it is easy
to see that the net $\{x^{*}_{\alpha}\}_{\alpha \in I}$ converges to
$j^{*}$ in the weak topology. Hence $x$ is a very smooth point of $X$. \hfill \ab
\end{proof}

\setcounter{rema}{4}
\begin{rema}
{\rm We have used the assumption $X/J$ has no smooth points to get a $j^{*}
\in J_1^{*}$. Thus this hypothesis can also be replaced by $d(x,J) <
1$. 
}
\end{rema}

\begin{rema}
{\rm If $J$ has property $(**)$ of [S] then arguments similar to the ones
given above can be used to show that under the same hypothesis every
smooth point of $X$ is a Fr\'{e}chet smooth point. 
}
\end{rema}

\setcounter{exam}{6}
\begin{exam}
{\rm 
Let $ 1 < p < \infty$ and $E$ be a Banach space such that $K (\ell^p,
E)$ is an $M$-ideal in ${\cal L} (\ell^p,E)$ (see Corollary 6.4 of
[KW] for a necessary and sufficient condition for this to happen). If
$T \in {\cal L} (\ell^p,E)$ is a smooth point, it follows from Theorem 1
of [GY2] that $\hbox{d}(T, K(\ell^p, E)) <1$. Thus if $E$ further satisfies
the condition $K(\ell^p,E)$ is Hahn--Banach smooth, then any smooth point
of ${\cal L}(\ell^p,E)$ is very smooth. 
}
\end{exam}

\begin{exam}
{\rm 
Let $X$ be an $M$-ideal in its bidual. See Chapters III and VI of [HWW]
for several examples of such spaces from among function spaces and
spaces of operators. It follows from our Proposition~2 and Theorem~2 of
[R7] where we have proved that $X$ is an $M$-ideal under appropriate
canonical embeddings in all the duals of even order, that every smooth
point of $X$ is very smooth and continues to be a very smooth point of
all the duals of even order of $X$. Assume further $X^{{*}{*}} / X$
has no smooth points then since $X$ is Hahn--Banach smooth (see [HWW],
Corollary~III.2.15), we get that any smooth point of $X^{{*}{*}}$ is
very smooth. $c_0 \subset \ell^{\infty}$ and ${\cal K}(\ell ^2) \subset
{\cal L}(\ell^2)$ are well-known examples of this phenomenon (see [KY]).
It may be worth recalling here that any such $X$ (see [HWW], 
Theorem~III.4.6) can be renormed to have a strictly convex dual (and hence every
unit vector is a smooth point) and continues to be an $M$-ideal in the
bidual with respect to the new norm. Thus in this renorming, all the
unit vectors are very smooth and continue to be very smooth in all the
duals of $X$ of even order. 
}
\end{exam}

The next example is once again from operator theory.

\begin{exam}
{\rm Let $X, Y$ be separable Banach spaces in the classes $(M_p)$ and $(M_q)$
for some $1 < p,q < \infty$ (see [HWW], Chapter VI and [W1]). Since
${\cal K}(X,Y)$ is an $M$-ideal in its bidual, ${\cal L}(X,Y)$ and since
the quotient space has no smooth points (see [HWW], Lemma~VI.5.18), it
follows that every smooth point of ${\cal L}(X,Y)$ is very smooth. In
particular in ${\cal L}(\ell^p)$ every smooth point is very smooth. It
follows from Theorem~6 of [R7] that there are no very smooth points in
$\ell^1$. It can thus be deduced from Proposition~17 in this paper that ${\cal
L}(\ell^1)$ does not have very smooth points. The question of very
smooth points of ${\cal L}(\ell^{\infty})$ will also be considered later
in this paper.
} 
\end{exam}

Our next example shows that unlike the situation described in the
previous examples, the fourth dual of $X$ (denoted by $X^{(IV)} $)
quotiented by $X$ can have smooth points. 

\begin{exam}
{\rm Let $X = c_0$. It follows from Theorem 2 of [R7] that it is an $M$-ideal
in $X^{(IV)}$. It is well known that $X^{(IV)}$ can be identified with a
$C(K)$ space and isolated points of $K$ correspond precisely to 
one-dimensional $L$-summands (or atoms) of $X^{{*}{*}{*}}$. Now $X = \{f
\in C(K) : f(K') = 0 \}$, where $K' = \{k \in K : f(k) = 0\ \hbox{for all}\  f
\in c_0 \}$. See Chapter~I of [HWW]. Thus $X^{(IV)} / X$ can be identified
as $C(K')$. To exhibit smooth points it is therefore enough to exhibit
isolated points in $K'$. Now for any $ \tau \in \partial_e (c_0^\bot)_1$
since $\hbox{span}\{\tau\}$ is an $L$-summand, it is an isolated point of $K
\cap K'$. 
}
\end{exam}

We do not know if a similar result is true in the non-commutative
situation of ${\cal K}(\ell^2) \subset {\cal L}(\ell^2)$. 

In the next proposition we again impose a condition on the quotient
space to preserve very smooth points and this result is modeled on
Proposition~3.4 of [GI]. Here we will be assuming that the subspace $J$
is only an ideal (in the sense of [GKS]) but make up for the loss in the
geometry of the projection in the dual, by making the topological
assumption that the quotient space is a Grothendieck space, that is,
weak$^{*}$ and weak sequential convergence coincide in the dual of the
quotient space. It is known that any von Neumann algebra is a
Grothendieck space, by [P].
We recall from [S] and [GI] that for a smooth point to be very smooth
one only needs to verify weak$^{*}$-weak sequential continuity of the
norm attaining functional. We also recall that any smooth Grothendieck
space is reflexive.

\setcounter{propo}{10}
\begin{propo}$\left.\right.$\vspace{.3pc}

\noindent Let $X$ be a smooth Banach space and let $J \subset X$ be an ideal such
that $X / J$ is a Grothendieck space. Then a very smooth point of $J$ is
also a very smooth point of $X$. 
\end{propo}

\begin{proof}
Let $P$ be a projection of norm one in $X^{*}$ such that
$\hbox{ker}\ P = J^\bot$. Let $j_0$ be a very smooth point of $J$ and let
$x_0^{*}$ be the unique functional with $x_0^{*}(j_0) = 1$. We shall
show that for any sequence $\{x^{*}_n\}_{n \geq 1} \subset X^{*}_1$
such that $x_n^{*} \rightarrow x_0^{*}$ in the weak$^{*}$-topology
also converges in the weak topology.

Since for any $x^{*}$, $P(x^{*})$ is the norm preserving extension
of $ x^{*} / J$, we have that $P(x_0^{*}) = x_0^{*}$ and also
$P(X^{*})$ is isometric to $J^{*}$ via the map $P(x^{*}) \rightarrow
x^{*} / J$. Thus $P(x_n^{*}) \rightarrow x_0^{*}$ in the 
weak$^{*}$-topology of $J^{*}$ and hence by hypothesis it would converge in the
weak topology. Therefore $\{P(x_n^{*}) - x_n^{*}\}_{n \geq 1} \subset
J^\bot$ converges to zero in the weak$^{*}$-topology and as $X / J$ is
a Grothendieck space, this convergence is also in the weak topology.
Hence $x_n^{*} \rightarrow x_0^{*}$ in the weak topology. Therefore,
$j_0$ is a very smooth point of $X$. \hfill \ab
\end{proof} 

We now recall a well-known characterization of smooth points of $K
(X,Y)$ due to Heinrich~[H].

\begin{thoe}[{\bf H}]. 
{\it \  $T \in K (X,Y)$ is a smooth point if and only if
$T^{*}$ attains its norm at a unique (up to a constant multiple)
$y_0^{*} \in \partial_ eY_1^{*}$ and $T^{*} y_0^{*}$ is a smooth
point of $X^{*}$. 
}\vspace{.5pc}
\end{thoe}

We next consider very smooth points of $K(X,Y)$. Let $T \in K(X,Y)$
be a smooth point. Note that since $T^{*}(y_0^{*})$ is a smooth point,
there exists a unique $x_0^{{*} {*}} \in \partial_ e X^{{*} {*}}_1$
such that $x_0^{{*} {*}} (T^{*} (y_0^{*}))=1$. Thus the linear
functional $x_0^{{*} {*}} \otimes y_0^{*}$ defined by $(x_0^{{*}
{*}} \otimes y_0^{*})(S)=x_0^{{*} {*}} (S^{*}(y_0^{*}))$ is an
extreme point of $K(X,Y)_1^{*}$ and is the unique functional attaining
its norm at $T$. Keeping the criterion from [GI] mentioned earlier in
view, we first look at points of $w^{*}-w$ continuity. As a further
geometric motivation for studying this concept, we recall from [HWW, 
p.~125] that $x^{*} \in X^{*}_1$ is a point of $w^{*}-w$ continuity if
and only if under the canonical embeddings, $x ^{*}$ has a unique norm
preserving extension to $X^{{*}{*}}$. 

\setcounter{theor}{11}
\begin{theor}[\!]
Let $x_0^{{*}{*}} \in \partial_e X^{{*}{*}}_1$ and let $y_0^{*} \in
\partial_e Y^{*}_1$. Suppose $x_0^{{*} {*}} \otimes y_0^{*} \in
\partial_ e K(X,Y)^{*}_1$ is a point of $w^{*} - w$ continuity. Then
$x_0^{{*} {*}}$ and $y_0^{*}$ are points of $w^{*} - w$ continuity.
In particular $x_0^{{*} {*}} = x_0 \in \partial_ e X_1$. 
\end{theor}

\begin{proof}
Let $\{x^{{*}{*}}_\alpha\}_{\alpha \in I}$ be a net in
$X_1^{{*} {*}}$ such that $x^{{*}{*}}_\alpha \stackrel
{w^{*}}{\rightarrow} x_0^{{*} {*}}$. For any $T \in K(X,Y),
(x^{{*}{*}}_\alpha \otimes y_0^{*})(T)= x^{{*} {*}}_\alpha
(T^{*}(y_0^{*})) \rightarrow x_0^{{*} {*}} (T^{*} (y_0^{*}))$.
Thus $x^{{*}{*}}_\alpha \otimes y_0 ^{*}
\stackrel{w^{*}}{\rightarrow} x_0^{{*} {*}} \otimes y_0^{*}$. Hence
by the hypothesis, $x^{{*}{*}}_\alpha \otimes y_0^{*}
\stackrel{w}{\rightarrow} x_0^{{*} {*}} \otimes y_0^{*}$. 

Now let $\tau \in X_1^{{*} {*} {*}}$ and let $y_0 \in Y$ be such that
$y_0^{*} (y_0)=1$. Consider the functional, $(\tau \otimes y_0)
(x^{{*}{*} } \otimes y^{*} ) = \tau(x^{{*}{*}}) y^{*} (y_0)$. It
follows from the well-known inclusions, $X^{{*}{*}{*}}
\otimes_{\epsilon} Y \subset (X^{*} \otimes_{\epsilon} Y)^{{*}{*}}
\subset {\cal K}(X, Y)^{{*}{*}}$, that $\tau \otimes y_0 \in K
(X,Y)^{{*} {*}}.$ Hence $(\tau \otimes y_0)(x^{{*}{*}}_\alpha
\otimes y^{*}_0) = \tau (x^{{*}{*}}_\alpha) \rightarrow \tau
(x_0^{{*} {*}})$. Hence $x^{{*}{*}}_\alpha \stackrel
{w}{\rightarrow} x_0^{{*} {*}}.$ Since $X_1$ is $w^{*}$-dense in
$X_1^{{*} {*}}$, we get that $x_0^{{*} {*}}=x_0 \in \partial_ e
X_1$. Similar arguments work for $y_0^{*}.$ 
\hfill \ab
\end{proof}

\setcounter{rema}{12}
\begin{rema}
{\rm It is clear from the proof that the arguments also go through for the
space ${\cal L}(X,Y)$. We note that for $x \in \partial_e
X^{{*}{*}}_1$ and $y^{*} \in \partial_e Y^{*}_1$, $x \otimes y^{*}$
in general need not be an extreme point of ${\cal L}(X,Y)^{*}_1$ (see
[HWW], p 267 for an example). 
}
\end{rema}

\setcounter{coro}{13}
\begin{coro}$\left.\right.$\vspace{.3pc}

\noindent Every very smooth point $T$ of $K(X,Y)$ attains its norm at a unique
unit vector {\rm (}up to constants{\rm )} $x_0$ and $T(x_0)$ is a very smooth point
of $Y$. 
\end{coro}

\begin{proof}
Let $T$ be a very smooth point of $K(X,Y)$. As noted
before, there exists a $x_0^{{*} {*}} \in \partial_ e X_1^{{*} {*}}$
and a $y_0^{*} \in \partial_ e Y_1^{*}$ such that $x_0^{{*} {*}}
\otimes y_0^{*}$ is the unique element of $\partial_ e K(X,Y)_1^{*}$
that attains its norm at $T$. Since $T$ is very smooth, $x_0^{{*} {*}}
\otimes y_0^{*}$ is a point of $w^{*}-w$ continuity. Therefore,
$x_0^{{*} {*}}=x_0 \in X$. Hence $1=y_0^{*} (T(x_0)) \leq \|T (x_0)
\| \leq \| T \|=1$. Thus $\| T \|=1=\|T(x_0) \|.$ That $T(x_0)$ is very
smooth again follows from the above theorem. \hfill \ab
\end{proof}

We do not know if the converse of the above theorem is always true.
However when $Y = C(K)$, for a compact set $K$ we have the following
complete description of very smooth points of ${\cal K}(X, C(K))$. It is
well known that ${\cal K}(X, C(K))$ can be idetified with the space of
vector-valued functions $C(K, X^{*})$. Since in this case the arguments
given during the proof of Theorem 12 are much simpler we present the
complete proof. 

\begin{coro}$\left.\right.$\vspace{.3pc}

\noindent Let $K$ be a compact Hausdorff space and let $X$ be any Banach space. $f
\in C(K,X)$ is very smooth if and only if there exists a unique point $k
\in K$ which is an isolated point such that $f(k)$ is a very smooth
point of $X$. 
\end{coro}

\begin{proof}
We begin by noting that the dual of $C(K,X)$ can be
identified with $M(K,X^{*})$, the space of $X^{*}$-valued regular
Borel measures equipped with the total variation norm.

Suppose $f \in C(K,X)$ is very smooth. Let $k \in K$ and $x^{*} \in
\partial_ e X_1^{*}$ be such that $1=x^{*} (f(k)).$ 

Let $\{x_\alpha^{*}\}_{\alpha \in I}$ be a net in $X_1^{*}$ such that
$x_\alpha^{*} \stackrel{w^{*}}{\rightarrow} x^{*}$. Consider for any
$\tau \in X^{{*} {*}}, F \in M(K, X^{*}), (\delta (k) \otimes
\tau)(F) =\tau(F(\{k\})).$ 
Then $\delta (k) \otimes \tau \in C(K,X)^{{*} {*}}$. Also $\delta (k)
\otimes x_\alpha^{*} \stackrel{w^{*}}{\rightarrow} \delta (k) \otimes
x^{*}$ and hence $\delta(k) \otimes x_\alpha^{*}
\stackrel{w}{\rightarrow} \delta(k) \otimes x^{*}.$ Therefore $(\delta
(k) \otimes \tau) (\delta (k) \otimes x_\alpha^{*})=\tau
(x_\alpha^{*}) \rightarrow (\delta (k) \otimes \tau)(\delta (k) \otimes
x^{*}) =\tau(x^{*}).$ Hence $x_\alpha^{*} \stackrel{w}{\rightarrow}
x^{*}.$ 

Also if $\{k_\alpha\}_{\alpha \in I}$ is a net in $K$ such that
$k_\alpha \rightarrow k$, then similar arguments show that $\delta
(k_\alpha) \stackrel{w}{\rightarrow} \delta (k)$ and hence $\{k\}$ is an
isolated point of $K$. Similar arguments show that $f(k)$ is a very
smooth point.

Conversely suppose that $k \in K$ is an isolated point such that $f(k)$
is a very smooth point. Since $M=\{ f \in C(K,X) : f(k) = 0\}$ is an
$M$-summand in $C(K,X)$, with $X$ (under the canonical identification)
as the complementary summand, it follows from Proposition 2 that $f$ is
a very smooth point.   \hfill \ab
\end{proof}

\begin{coro}$\left.\right.$\vspace{.3pc}

\noindent Any very smooth point of ${\cal K} (X,C(K))$ is a very smooth point of
${\cal L}(X,C(K))$.
\end{coro}

\begin{proof}Since $C(K)$ has the MAP we have from our earlier
remark that a very smooth point of ${\cal K} (X,C(K))$ is a smooth point
of ${\cal L}(X,C(K))$. We now use the identification of ${\cal L}(X,
C(K))$ as the space $W^{*} C(K, X ^{*})$ of functions on $K$, that are
continuous when $X^{*}$ is equipped with the $w^{*}$-topology,
equipped with the supremum norm. Suppose $f \in C(K,X^{*})$ is very
smooth. There exists an isolated point $k_0 \in K$ and an $x_0 \in
\partial_ e X_1$ that is a point of $w^{*} - w$ continuity of
$X_1^{{*} {*}}$ such that $f(k_0)(x_0)=1$. Since $f \rightarrow f
\chi_{\{k_0\}}$ is an $M$-projection in $W^{*} C(K,X^{*})$, once more
using arguments similar to the ones given during the proof of the
previous corollary, we get that $f$ is a very smooth point of $W^{*}
C(K,X^{*}).$ \hfill \ab
\end{proof}

It follows from the arguments given above that if a very smooth point
$f \in W^{*} C(K, X^{*})$ attains its norm at a $k \in K$ and $f(k)$
attains its norm, then $\{ k\}$ is an isolated point and $f(k)$ is a
very smooth point of $X^{*}$. We use similar ideas in the next
proposition to describe very smooth points in $\ell^{\infty}$-direct
sums of Banach spaces. For any set $I$, if $\beta(I)$ denotes the
Stone--Cech compactification of the discrete space $I$, the space $W^{*}
C(K, X^{*})$ can be identified as the $\ell^{\infty}$ direct sum of
$|I|$-many copies of $X^{*}$. Our arguments run parallel to those given
during the proof of Theorem 1 in [GY1]. 

\setcounter{propo}{16}
\begin{propo}$\left.\right.$\vspace{.3pc}

\noindent Let $\{X_i\}_{ i \in I}$ be a family of Banach spaces. Let $X =
\oplus_{\infty} X_i$.~ $x \in X$ is a very smooth point iff there exists
a unique $i_0$ such that $\|x(i_0)\| = 1\! >\!\hbox{\rm Sup}\ \{ \|x(i)\| : i \neq i_0
\}$ and $x(i_0)$ is a very smooth point of $X_{i_0}$. 
\end{propo}

\begin{proof}
Suppose there is no $i_0$ such that $\|x(i_0)\| = 1 $ or
there is an $i_0$ such that $\|x(i_0)\| = 1 =\ \hbox{Sup}\ \{ \|x(i)\| : i \neq
i_0 \}$. In either case arguing as in the proof of Theorem~1 in [GY1],
we get a disjoint partition $\{A,B\}$ of $I$, such that the supremum of
$\|x(i)\|$ over both sets is $1$. Now using the canonical projection
into the summands of $X$ formed out of $X_i$'s taken from $A$ and $B$,
we get a decomposition $x = y + z$ with $\|y \| = 1 = \|z \|$ and such
that $\|y - z \| \leq 1$. Thus it follows from Lemma 2 of [GY1] that $x$
is not a smooth point. This contradiction establishes that there is a
unique $i_0$ such that $\|x(i_0)\| = 1$. Since $X_{i_0}$ is an
$M$-summand and hence an $M$-ideal of $X$, the conclusions follow from
Proposition~2. \hfil \ab
\end{proof}

Inherent in the discussion of this section is that a very smooth
point of $C(K)$ is the indicator function of an isolated point and hence
is a Fr\'{e}chet smooth point. We use this in the next proposition to
show that for a $JB^{*}$-triple (see [MR] for the definitions) very
smooth points and Fr\'{e}chet smooth points coincide. 

\begin{propo}$\left.\right.$\vspace{.3pc}

\noindent  Let $X$ be a $JB^{*}$-triple. Any very smooth point
of $X$ is a Fr\'{e}chet smooth point.
\end{propo}

\begin{proof}
Let $x \in X$ be a very smooth point. Let $C(x)$ be the
$JB^{*}$-subtriple generated by odd powers of $x$. By Proposition 2 we
get that $x$ is a very smooth point of $C(x)$. Since $C(x)$ is a
commutative $C^{*}$-algebra, we get that $x$ is a Fr\'{e}chet smooth
point. Since the sequence of odd powers of $x$ converge in the
weak$^{*}$-topology of the bidual to a tripotent $u(x)$ and as the norm
closed faces $\{x\}$ and $\{u(x)\}$ coincide (see [MR], Lemma~2.4) we
conclude as in the proof of Theorem~3.1 in [MR] that $x$ is a
Fr\'{e}chet smooth point of $X$. \hfill \ab
\end{proof}

We conclude this section with a renorming result for very smooth
points. We have from the results in \S3 of [S] that if $X$ is very
smooth then $X^{*}$ has the Radon--Nikod\'{y}m property and if $X^{*}$
is very smooth then $X$ is reflexive. Also if $X$ has the
Radon--Nikod\'{y}m property then $X^{*}$ has a point of Fr\'{e}chet
differentiability and hence a very smooth point (see [B], 
Theorem~5.7.4). It is therefore natural to ask for renormings that admit very
smooth points. 

\begin{propo}$\left.\right.$\vspace{.3pc}

\noindent  Suppose for every renorming of $X$ the dual unit
ball has a very smooth point. Then $X$ has the Radon--Nikod\'{y}m
property.
\end{propo}

\begin{proof}
Suppose $x^{*} \in X^{*}_1$ is a very smooth point. As
noted before, there exists a unit vector $x$ such that $x^{*}(x) = 1 =
\|x^{*}\|$. We claim that $x$ is an extreme point of the unit ball of
the fourth dual $X^{(IV)}$ of $X$. Suppose $x = \frac{1}{2} \{\lambda_1
+ \lambda_2\}$ where $\lambda_i \in X^{(IV)}_1$. Thus $1 = x^{*}(x) =
\lambda_1 (x^{*}) = \lambda_2(x^{*})$. Since $x^{*}$ is a smooth
point of $X^{{*}{*}{*}}$ we get that $x = \lambda_1 = \lambda_2$.
Therefore $x \in \partial_e X^{(IV)}_1$. Thus for every equivalent norm
on $X$, the unit ball has an extreme point of the unit ball (w. r. t the
equivalent norm) of $X^{(IV)}$. Hence it follows from Corollary 6 of
[Hu] that $X$ has the Radon--Nikod\'{y}m property. \hfill \ab
\end{proof}

\section{Connection with the work of Ruess and Stegall} 

In a series of papers in the 80's Ruess and Stegall ([RS1], [RS2]) have
showed that $w^{*}$-denting points of ${\cal
K}(X,Y)^{*}_1$ and ${\cal L}(X,Y)^{*}_1$ coincide and are precisely
points of the form $x^{{*} {*}} \otimes y^{*}$, where $x^{{*} {*}}$
and $y^{*}$ are $w^{*}$-denting points of $X^{{*} {*}}_1$
and $Y_1^{*}$ respectively. Using a result of Lin~{\it et~al} [LLT], 
a denting ($w^{*}$-denting) point is precisely an extreme point
that is a point of weak-norm ($w^{*}$-norm) continuity of the identity
map in the dual unit ball. 

~~While investigating the question of very smooth points in the unit
ball of spaces of operators one encounters extreme points of ${\cal
L}(X,Y)_1^{*}$ (or $K(X,Y)_1^{*}$) that are points of $w^{*}-w$
continuity. Our theorem in the previous section is an attempt to prove a
Ruess--Stegall type result in one direction (see Remark 13). We could
describe them fully only for the space of compact operators when $Y
=C(K)$. Next proposition gives another partial answer and should again
be compared with Remark~13. 

\begin{propo}$\left.\right.$\vspace{.3pc}

\noindent  Suppose $X^{*}$ or $Y$ has the {\rm CMAP}. Let $x \otimes
y^{*} \in \partial_e {\cal K}(X,Y)^{*}_1$ be a $w^{*}-w$ point of
continuity. Then $x \otimes y^{*} \in \partial_e {\cal L}(X,Y)^{*}_1$.
\end{propo} 

\begin{proof}
 As noted before the hypothesis implies that ${\cal
K}(X,Y) \subset {\cal L}(X,Y) \subset {\cal K}(X,Y)^{{*}{*}}$.
Applying Lemma~III.2.14 from [HWW] once more we have that $x \otimes
y^{*} $ has a unique norm preserving extension to a functional on
${\cal K}(X,Y)^{{*}{*}}$. Since $x \otimes y^{*} \in {\cal
L}(X,Y)^{*}$ is one such extension and as $x \otimes y^{*} \in
\partial_e {\cal K}(X,Y)^{*}_1$, we get that $x \otimes y^{*} \in
\partial_e {\cal L}(X,Y)^{*}_1$. \hfill\ab
\end{proof}

In the next proposition we show that under the $M$-ideal assumption
the dual unit balls of the space of both compact and bounded operators
have the same $w^{*}-w$ points of continuity. 

\begin{propo}$\left.\right.$\vspace{.3pc}

\noindent  Let ${\cal K}(X,Y) \subset {\cal L}(X,Y)$ be an
$M$-ideal. Any $\tau \in \partial_e {\cal L}(X,Y)^{*}_1$ that is a
$w^{*}-w$ point of continuity is of the form $x \otimes y^{*}${\rm ,} where
$x \in \partial_e X^{{*}{*}}_1$ is a $w^{*}-w$ point of continuity
and $y^{*} \in \partial_e Y^{*}_1$ is also a $w^{*}-w$ point of
continuity.
\end{propo} 

\begin{proof}
Let $\tau \in \partial_e {\cal L}(X,Y)^{*}_1$ be a
$w^{*}-w$ point of continuity. We shall show that $\tau \in \partial_e
{\cal K}(X,Y)^{*}_1$. It would then follow from the arguments given
during the proof of Lemma 1 in [R5] that it is a $w^{*}-w$ point of
continuity of ${\cal K}(X,Y)^{*}_1$. Hence the conclusion follows from
Theorem~12.

It follows from the arguments given during the proof of Proposition 1
in [R4] that a net of convex combinations of functionals of the form
$\{x \otimes y^{*} : \|x\| = 1 = \|y^{*}\| \}$ ( which are in ${\cal
K}(X,Y)^{*}_1 \subset {\cal L}(X,Y)^{*}_1)$ converges to $\tau$ in the
weak$^{*}$ topology of ${\cal L}(X,Y)^{*}$. Hence by our assumption
this convergence occurs in the weak topology. Therefore $\tau \in
\partial_e {\cal K}(X,Y)^{*}_1$. \hfill \ab 
\end{proof}

The following corollary is now easy to see and should be compared
with the description of smooth points where one needed extra assumptions
on the quotient space, see [KY], [GY2] and [HWW] Chapter VI. 

\setcounter{coro}{21}
\begin{coro}$\left.\right.$\vspace{.3pc}

\noindent  Let ${\cal K}(X,Y) \subset {\cal L}(X,Y)$ be an
$M$-ideal. If $T \in {\cal L}(X,Y)$ is a very smooth point then it
attains its norm at a unique vector $x$ and $T(x)$ is a very smooth
point of $Y$. 
\end{coro}

\section{Spaces without  points of $\hbox{\sa w}^{*}-\hbox{\sa w}$ continuity}

In the concluding part of this paper we consider situations where
$X^{*}_1$ fails to have points of $w^{*}-w$ continuity. For such
spaces, $X_1$ fails to have very smooth points. We recall from [HWW]
that $Y \subset Z$ is said to be an $L$-ideal if there is an onto
projection $P: Z \rightarrow Y$ such that $\|z\| = \|P(z)\| + \|z -
P(z)\|$ for all $z \in Z$. Also if a Banach space $X$, when considered
under the canonical embedding, as a subspace of $X^{{*}{*}}$, is an
$L$-ideal then $X$ is said to be an $L$-embedded space. See [HWW],
Chapter~IV and VI for several examples of $L$-embedded spaces from
among function spaces and spaces of operators. In what follows we only
consider non-reflexive spaces. Our first result leads to an example of a
smooth space that does not have any very smooth points. 

\setcounter{propo}{22}
\begin{propo}$\left.\right.$\vspace{.3pc}

\noindent  Let $X \subset Y \subset X^{{*}{*}}$ {\rm (}under the
canonical embedding{\rm )} and suppose $X$ is a proper $L$-ideal in $Y$. Then
$X^{*}_1$ has no $w^{*}-w $ points of continuity. In particular if $X$
is an $L$-embedded space then $X_1$ has no very smooth points.
\end{propo} 

\begin{proof}
Let $P$ be the $L$-projection in $Y$ whose range is $X$.
Let $x^{*} \in X^{*}_1$ be a $w^{{*} }-wPC$. We get the necessary
contradiction by showing that $x^{*}$ has no unique norm preserving
extension to $Y^{*}$ and hence to $X^{{*}{*}}$. Choose a $\tau \in
X^{\bot}$, i.e., such that $P^{{*} }(\tau ) =0$ and $ 0<\Vert \tau \Vert
\leq \Vert P^{{*} }(x^{*})\Vert \leq \Vert x^{*}\Vert $. We note that
since $X \subset Y \subset X^{{*}{*}}$ under the canonical embedding,
$x^{*}$ has a natural extension as the evaluation functional to $Y$. We
continue to denote this extension by $x^{*}$. Now for any $x \in X$,
$(P^{{*} }(x^{*})+\tau )(x) = x^{{*} }(x)$ and as $P^{{*} }$ is a
$M$-projection, $\Vert P^{{*} }(x^{*})+\tau \Vert =\max\{\Vert P^{{*}
}(x^{*})\Vert ,\Vert \tau \Vert \}\leq \Vert x^{*}\Vert $. Hence there
is no unique extension. \hfill \ab
\end{proof}

\setcounter{rema}{23}
\begin{rema}
{\rm It is well known that the Hardy space on the unit circle,
$H^1 _0$, is a smooth $L$-embedded space (see [HWW], p. 167). However,
the unit ball has no very smooth points.
}
\end{rema}

In view of this remark the following question can be asked: `does
there exist a very smooth space $X$ so that no unit vector is a smooth
point of the fourth dual of $X$?' This however is not the case. To see
this note that since $X$ is very smooth, $X^{*}$ has the
Radon--Nikod\'{y}m property (see [S], \S3), thus $X$ has points of
Fr\'{e}chet differentiability (see [B], Proposition~5.6.13). It is well
known that a point of Fr\'{e}chet differentiability continues to be a
point of Fr\'{e}chet differentiability of all the duals of even order of
$X$. 

\setcounter{coro}{24}
\begin{coro}$\left.\right.$\vspace{.3pc}

\noindent Let $X$ be an $L$-embedded space with the metric
approximation property and let $Y$ be any Banach space. There are no
points of $w^{*}-w$ continuity in the dual unit ball of the projective
tensor product space $X \otimes_{\pi} Y$.
\end{coro}

\begin{proof}
We use the known inclusions (under canonical embeddings)
\begin{equation*}
X \otimes_{\pi} Y \subset X^{{*}{*}} \otimes_{\pi} Y \subset (X
\otimes_{\pi} Y)^{{*}{*}}.
\end{equation*}
See [R6]. Since $X$ is an $L$-ideal in
$X^{{*}{*}}$ it follows from Theorem VI.6.8 of [HWW] that $X
\otimes_{\pi} Y$ is an $L$-ideal in $X^{{*}{*}} \otimes_{\pi} Y$. It
thus follows from our above proposition that there are no points of
$w^{*}-w$ continuity in the dual unit ball. \hfill \ab
\end{proof}

\setcounter{rema}{25}
\begin{rema}
{\rm For any positive measure $\mu$, $L^1(\mu)$ is an
$L$-embedded space with the metric approximation property. For any
Banach space $Y$, let $L^1(\mu,Y)$ denote the space of $Y$-valued
Bochner integrable functions. Since $L^1(\mu) \otimes_{\pi} Y =
L^1(\mu,Y)$ we get that there are no $w^{*}-w$ points of continuity in
the dual unit ball of $L^1(\mu,Y)$. We recall that ${\cal L}(X,Y^{*})$
can be identified with $(X \otimes_{\pi} Y)^{*}$. Thus when $X$ is an
$L$-embedded space with the metric approximation property there are no
points of $w^{*}-w$ continuity in ${\cal L}(X,Y^{*})_1$.
}
\end{rema} 

Let $X$ be a Banach space. Suppose there is a $u \in X^{{*}{*}}
\backslash X $ such that $\|u + x\| = \|u\| + \|x\|$ for every $x \in
X$. Then $Y = \hbox{span}\ \{X,u\}$ satisfies the hypothesis of the above
proposition. This geometric condition has been well studied in the
literature, see [DGZ] . We use these ideas in the next corollary. 

Turning once more to spaces of operators we recall that a Banach
space $X$ is said to have the Daugavet property if $\|I + T\| = 1 +
\|T\|$ for all compact operators $T$ (see [W], Definition~2.1 and Theorem
2.7). Any space with the Daugavet property contains an isomorphic copy
of $\ell^1$ (see [W], Theorem~2.6). 

\setcounter{coro}{26}
\begin{coro}$\left.\right.$\vspace{.3pc}

\noindent  Let $X$ be a Banach space having the Daugavet
property and the metric approximation property. Then ${\cal
K}(X)^{*}_1$ has no points of $w^{*}-w $ continuity.
\end{coro}

\begin{proof}
Since $X$ has the metric approximation property it follows
from Example 1 in [R6] that ${\cal K}(X) \subset {\cal L}(X) \subset
{\cal K} (X)^{{*}{*}}$ under the canonical embedding. Since $X$
satisfies the Daugavet property the conclusion follows. \hfill \ab
\end{proof}

This author has proved recently ([R2], [R3]) that when $K$ is
infinite, for the space ${\cal L}(X,C(K))$ there are no point of
$w$-norm continuity in the unit ball. When ${\cal L}(X,C(K))$ is a dual
space it would be interesting to know if the unit ball can have points
of $w^{*}-w$ continuity (with respect to a predual!) . We have some
partial results. 

\setcounter{propo}{27}
\begin{propo}$\left.\right.$\vspace{.3pc}

\noindent  Let $X$ be a reflexive Banach space. If ${\cal L}(X,
C(K))$ is a dual space then the unit ball has no points of $w^{*}-w$
continuity. In particular the predual has no very smooth points.
\end{propo}

\begin{proof} It follows from [CG] that when $X$ is reflexive the
assumption ${\cal L}(X, C(K))$ is a dual space implies that $K$ is
hyperstonean. Hence the (unique) predual is of the form $L^1
(\vartheta,X)$ for a category measure $\vartheta$ on $K$. Since $X$ is
reflexive, it follows from p. 200 of [HWW], that $L^1(\vartheta,X)$ is an
$L$-embedded space. Hence ${\cal L} (X,C(K))_1$ has no points of
$w^{*}-w$ continuity. \hfill \ab
\end{proof}

\setcounter{rema}{28}
\begin{rema}
{\rm Similar arguments work when $X$ is the predual of a von
Neumann algebra. In this case if ${\cal L}(X, C(K))$ is a dual space,
then $K$ is extremally disconnected and hence ${\cal L}(X, C(K))$ is a
von Neumann algebra and hence is the dual of an $L$-embedded space.
}
\end{rema} 

\begin{rema}
{\rm Turning to the case of weakly compact operators valued in
a $C(K)$ space or more generally if one considers the space $WC(K, X)$,
(functions continuous when $X$ has the weak topology), this author has
proved in [R1] that if $WC(K,X)$ is a dual space then $X$ is reflexive
and $K$ is hyperstonean. Therefore, it follows from the above proposition
that when $WC (K, X)$ is a dual space the unit ball has no $w^{*}-w$
PC's.
}
\end{rema}

\end{document}